 \documentclass[final,3p,times]{elsarticle}

\usepackage{amssymb}
 \usepackage{amsthm}
\usepackage{amsmath,amssymb,amsopn,amsfonts,mathrsfs,amsbsy,amscd}

\newcommand{\nil}{\mathrm{Nil} }
\newcommand{\sol}{\mathrm{Sol} }

\newcommand{\prs}{\langle\;,\;\rangle}

\newcommand{\too}{\longrightarrow}

\newcommand{\esp}{\quad\mbox{and}\quad}
\newcommand{\wi}{\widetilde}

\newcommand{\G}{\mathfrak{g}}
\newcommand{\g}{\mathfrak{g}}

\newcommand{\tr}{{\mathrm{tr}}}
\newcommand{\ric}{{\mathrm{ric}}}
\newcommand{\Ric}{{\mathrm{Ric}}}

\newcommand{\B}{{\cal B}}

\newcommand{\al}{\alpha}
\newcommand{\be}{\beta}

\newcommand{\e}{\epsilon}

\newcommand{\la}{\lambda}

\font\bb=msbm10

\def\B{\hbox{\bb B}}
\def\R{\hbox{\bb R}}

\newtheorem{theo}{Theorem}[section]
\newtheorem{pr}{Proposition}[section]

\newtheorem{rem}{Remark}

\begin{document}

\begin{frontmatter}


 
 \title{  Lorentzian left invariant metrics on three dimensional unimodular Lie groups and their curvatures}
 
 \author[label1,label2]{ Mohamed Boucetta, Abdelmounaim Chakkar}
 \address[label1]{Universit\'e Cadi-Ayyad\\ Facult\'e des sciences et techniques\\
 	BP 549 Marrakech Maroc.\\e-mail:  m.boucetta@uca.ac.ma  }
 \address[label2]{Universit\'e Cadi-Ayyad\\
 	Facult\'e des sciences et techniques\\
 	BP 549 Marrakech Maroc\\e-mail: abdelmounaim.chakkar@edu.uca.ma
 }


\begin{abstract}There are five unimodular  simply connected three dimensional unimodular non abelian Lie groups:  the nilpotent Lie group $\nil$, the special unitary group $\mathrm{SU}(2)$, the universal covering group $\wi{\mathrm{PSL}}(2,\R)$ of the special linear group, the solvable Lie group $\sol$ and the universal covering group $\wi{\mathrm{E}_0}(2)$ of the connected component of the Euclidean group. For each $G$ among these Lie groups, we give explicitly the list of all Lorentzian left invariant metrics on $G$, up to un automorphism of $G$. Moreover, for any Lorentzian left invariant metric in this list we give its Ricci curvature, scalar curvature, the signature of the Ricci curvature and we exhibit some special features of these curvatures. Namely, we give all the metrics with constant curvature, semi-symmetric non locally symmetric metrics and the Ricci solitons.  
\end{abstract}


\end{frontmatter}






\section{Introduction}\label{section1}
This paper has two goals:\begin{enumerate}
	\item  to complete the study done by Rahmani in \cite{rahmani} by classifying for each three dimensional simply connected unimodular Lie group $G$ all the Lorentzian left invariant  metrics on $G$, up to automorphism of $G$ and hence achieve a similar study to the one  done in the Riemannian case in \cite{lee},
	\item to study, for each class of these metrics, their curvatures.
\end{enumerate}
There are five unimodular non abelian  simply connected three dimensional unimodular Lie groups characterized by the signature of their Killing form:  the nilpotent Lie group $\nil$ with signature $(0,0,0)$, the special unitary group $\mathrm{SU}(2)$ with signature $(-,-,-)$, the universal covering group $\wi{\mathrm{PSL}}(2,\R)$ of the special linear group with signature $(+,+,-)$, the solvable Lie group $\sol$ with signature $(+,0,0)$ and the universal covering group $\wi{\mathrm{E}_0}(2)$ of the connected component of the Euclidean group with signature $(-,0,0)$. Let $(G,h)$ one of these Lie groups endowed with a Lorentzian left invariant  metric and $\g$ its Lie algebra with a fixed orientation. Denote by $\prs$ the values of $h$ at the neutral element. The study of Rahmani in \cite{rahmani} is based on a remark first made by Milnor in \cite{milnor}. There exists a product $\times:\g\times\g\too\g$ depending on $\prs$ and the orientation and $L:\g\too\g$ a symmetric endomorphism such that, for any $u,v\in\g$, the Lie bracket on $\g$ is given by
\begin{equation}\label{br}
[u,v]=L(u\times v).
\end{equation}Note that $L$ changes to $-L$ when we change the orientation. It is well-known (see \cite{oneil}) that there are four types of symmetric endomorphisms on a Lorentzian vector space. Depending on the type of $L$, there exists $\B_1=(e_1,e_2,e_3)$ an orthonormal basis of $\g$ with $\langle e_1,e_1\rangle=\langle e_2,e_2\rangle=1$ and $\langle e_3,e_3\rangle=-1$ such that \eqref{br} gives one of the following forms::
\begin{enumerate}
	\item Type $\mathrm{diag}\{a,b,c\}$. There exits $a,b,c\in\R$ such that
	\begin{equation}\label{bd}
	[e_1,e_2]=-ce_3, [e_2,e_3]=ae_1 \esp [e_3,e_1]=b e_2.
	\end{equation}
	In this case the eigenvalues of the matrix of the Killing form in $\B_1$  are $[-2ab,2ac,2bc]$.
	\item Type  $\{az\bar{z}\}$. There exits $a,\al,\be\in\R$ with $\be\not=0$ such that
	\begin{equation}\label{bazz}
	[e_1,e_2]=-\be e_2-\al e_3,\;[e_2,e_3]=ae_1\esp  [e_3,e_1]=\al e_2-\be e_3. 
	\end{equation}
	In this case the eigenvalues of the matrix of the Killing form in $\B_1$  are
	$[2(\al^2+\be^2),2a\sqrt{\al^2+\be^2},-2a\sqrt{\al^2+\be^2}]$.
	\item Type  $\{ab2\}$. There exits $a,b\in\R$   such that
	\begin{equation}\label{bab2}
	[e_1,e_2]=\frac12 e_2+\left(\frac12-b\right)e_3,\;[e_2,e_3]=ae_1\esp  [e_3,e_1]=\left(b+\frac12\right)e_2+\frac12e_3. 
	\end{equation}
	In this case the eigenvalues of the matrix of the Killing form in $\B_1$  are
	$[2b^2,(\sqrt{4b^2+1}-1)a,-(\sqrt{4b^2+1}+1)a]$.
	 \item Type  $\{a3\}$. There exits $a\in\R$   such that
	\begin{equation}\label{ba3}
	[e_1,e_2]=\frac1{\sqrt{2}}e_2-ae_3,\;[e_2,e_3]=ae_1+\frac1{\sqrt{2}}e_2\esp  [e_3,e_1]=\frac1{\sqrt{2}}e_1+ae_2+\frac1{\sqrt{2}}e_3. 
	\end{equation}In this case the eigenvalues of $M(K,\B_1)$ are complicated. But $\det M(K,\B_1)=-8a^6$ and $\tr(M(K,\B_1))=2(a^2+1)$ which can be used to determine the signature of $K$.
\end{enumerate} Moreover,  the type of $L$ is an invariant of Lorentzian left invariant metrics on unimodular three dimensional Lie groups. Indeed,
 let $G$ be an unimodular three dimensional Lie group and $h_1$ and $h_2$ are two left invariant Lorentzian metrics on $G$ and $L_1$ and $L_2$ the associated endomorphisms for a fixed orientation on $\g$. If there exists an automorphism $\phi:G\too G$ such that $h_1=\phi^*(h_2)$ then $L_2=\pm\psi\circ L_1\circ \psi^{-1}$ where $\psi=T_e\phi$. Thus if $L_1$ and $\pm L_2$ have two different types or have the same type and two different sets of eigenvalues then $h_1$ and $h_2$ are not isometric.

	Based on what above, our method in order to fulfill our  goals mentioned above goes as follows:
	\begin{enumerate}
		\item We consider a three unimodular Lie group $G$ and we fix a natural basis $\B_0$ of its Lie algebra $\g$ (see Section \ref{section2} where these bases are given).
		\item We endow $G$ with a left invariant Lorentzian metric $h_0$ and we denote by $L$ its associated endomorphism.
		\item Depending on the signature of the Killing form we can determine which are the possible types of $\pm L$.
		\item For each possible type there exists an orthonormal basis $\B_1=(e_1,e_2,e_3)$ such that the Lie bracket has one of the forms \eqref{bd}-\eqref{ba3}.
		\item  We find a basis $\B_2=(xe_1+ye_2+ze_3,ue_1+ve_2+we_3,pe_1+qe_2+re_3)$ of $\g$ which has the same constant structures as $\B_0$ and
		we consider the automorphism of Lie algebra $P:\g\too\g$ which sends the basis $\B_0$ to $\B_2$ and $\phi$  the automorphism of $G$ associated to $P$.  We put  $h_1=\phi^*(h_0)$.

		\item In some cases, we use another automorphism $\phi_1$ of $G$ such that $\phi_1^*(h_1)$ has a more reduced form than $h_1$.
		\item Finally, we give the matrix of $\phi_1^*(h_1)$ in $\B_0$. The $\phi_1^*(h_1)$ obtained constitute a list of Lorentzian left invariant metrics on $G$ depending on a reduced number of parameters and each Lorentzian left invariant metric on $G$ is isometric to one in this list. We find twenty non isometric classes of such metrics which shows that the situation is far more rich than the Riemannian case \cite{lee}.
		\item  We compute for each metric in the list the Ricci curvature and the scalar curvature which determine all the curvature since we are in dimension 3.  Finally,  we exhibit some metrics with distinguished curvature properties.
	\end{enumerate}
	
	The steps 5. and 6. are  the most difficult and, to find the basis $\B_2$ or the automorphism $\phi$, we have used the softwear Maple and the expression of the groups of automorphisms of unimodular three dimensional Lie algebras given in \cite{lee}. All the computations have been checked by Maple.

	 The paper is organized as follows. In Section \ref{section2},  we precise the models of unimodular three dimensional Lie groups we will use in this paper. In Section \ref{section3}, we perform for each unimodular three dimensional Lie groups the steps mentioned above and we give its list
	 of left invariant Lorentzian metrics. For $\nil$ the list contains three non isometric classes of metrics, for $\mathrm{SU}(2)$  one class, for $\wi{\mathrm{PSL}}(2,\R)$   seven non isometric classes, for $\sol$  seven non isometric classes and for
	 $\wi{\mathrm{E}_0}(2)$  three non isometric classes. These metrics are given by the formulas \eqref{n1}-\eqref{e2a02}.
	 In Section \ref{section4}, we give for each class of metric found in Section \ref{section3} its    Ricci curvature,  scalar curvature and the signature of the Ricci curvature. We give a table describing the possible signature of the Ricci curvature and the metrics realizing these signatures and we give the different types of Ricci operators (see Proposition \ref{pr}).
	 Finally, we recover some known results. Namely, Lorentzian left invariant metrics on unimodular three dimensional Lie groups which are of constant curvature, Einstein, locally symmetric, semi-symmetric not locally symmetric or Ricci soliton have been determined in \cite{boucetta, bro, calvaruso1} by giving their Lie algebras as in \eqref{bd}-\eqref{ba3}. We give  their corresponding metrics  in our list (see Theorems \ref{main1}-\ref{main3}).

\section{Preliminaries}\label{section2}
Let $(G,h)$ be a Lie group endowed with a Lorentzian left invariant metric $h$, $\g$ its Lie algebra and $\prs=h(e)$.
The
Levi-Civita connection of $(G,h)$ defines a product $\mathrm{L}:\G\times\G\too\G$ called the Levi-Civita product and given by  Koszul's
formula
\begin{eqnarray}\label{levicivita}2\langle
\mathrm{L}_uv,w\rangle&=&\langle[u,v],w\rangle+\langle[w,u],v\rangle+
\langle[w,v],u\rangle.\end{eqnarray}
For any $u,v\in\G$, $\mathrm{L}_{u}:\G\too\G$ is skew-symmetric and $[u,v]=\mathrm{L}_{u}v-\mathrm{L}_{v}u$. 
The curvature on $\G$ is given by
$
\label{curvature}\mathrm{K}(u,v)=\mathrm{L}_{[u,v]}-[\mathrm{L}_{u},\mathrm{L}_{v}].
$ The Ricci tensor is the symmetric tensor $\ric$ given by
$\ric(u,v)=\tr(w\too \mathrm{K}(u,w)v)$ and the Ricci operator $\Ric:\g\too\g$ is given by the relation $\langle \Ric(u),v\rangle=\ric(u,v)$. The scalar curvature is given by $\mathfrak{s}=\tr(\Ric)$. Recall that:
\begin{enumerate}\item $(G,h)$ is called flat if $K=0$;
	\item $(G,h)$ has constant sectional curvature if there exists a constant $\la$ such that, for any $u,v,w\in\g$,
	\[ \mathrm{K}(u,v)w=\la\left( \langle v,w\rangle u-\langle u,w\rangle v  \right). \]
\item $(G,h)$ is called Einstein if there there exists a constant $\la$ such that $\Ric=\la \mathrm{Id}_\g$.
\item $(G,h)$ is called locally symmetric if, for any $u,v,w\in\g$,
\[ \mathrm{L}_u(\mathrm{K})(v,w):=[\mathrm{L}_u,\mathrm{K}(v,w)]-\mathrm{K}(\mathrm{L}_uv,w)-\mathrm{K}(v,\mathrm{L}_uw)=0. \]
\item $(G,h)$ is called  semi-symmetric if, for any $u,v,a,b\in\g$,
\[ [\mathrm{K}(u, v),\mathrm{K}(a, b)]=\mathrm{K}(\mathrm{K}(u, v)a,b)+\mathrm{K}(a,\mathrm{K}(u,v)b). \]
\end{enumerate}

They are five unimodular  simply connected three dimensional unimodular non abelian Lie groups: \begin{enumerate}
	
	\item The nilpotent Lie group $\nil$ known as Heisenberg group whose Lie algebra will be denoted by $\mathfrak{n}$. We have
	\[ \nil=\left\{\left(\begin{matrix} 1&x&z\\0&1&y\\0&0&1\end{matrix}\right),x,y,z\in\R   \right\}
	\esp \mathfrak{n}=\left\{\left(\begin{matrix} 0&x&z\\0&0&y\\0&0&0\end{matrix}\right),x,y,z\in\R   \right\} \]and the non-vanishing Lie brackets in the canonical basis $\B_0=(X,Y,Z)$ are given by $[X,Y]=Z$. The Killing form is trivial.
	\item $\mathrm{SU}(2)=\left\{\left(\begin{matrix} a+bi&-c+di\\c+di&a-bi\end{matrix}\right),a^2+b^2+c^2+d^2=1   \right\}
	\esp \mathfrak{su}(2)=\left\{\left(\begin{matrix} iz&y+ix\\-y+xi&-zi\end{matrix}\right),x,y,z\in\R   \right\}.$ In the basis
	$$\B_0=\left\{\sigma_x=\left(\begin{matrix}0&i\\i&0\end{matrix} \right), \sigma_y=\left(\begin{matrix}0&1\\-1&0\end{matrix} \right), \sigma_z=\left(\begin{matrix}-i&0\\0&i\end{matrix} \right)\right\},$$ we have
	\begin{equation} \label{eqsu2}[\sigma_x,\sigma_y]=2\sigma_z,\;[\sigma_y,\sigma_z]=2\sigma_x\esp [\sigma_z,\sigma_x]=2\sigma_y. \end{equation} The Killing form has signature $(-,-,-)$.
	
	\item The universal covering group $\wi{\mathrm{PSL}}(2,\R)$ of $\mathrm{SL}(2,\R)$ whose Lie algebra  is $\mathrm{sl}(2,\R)$. In the basis
	$$\B_0=\left\{X_1=\left(\begin{matrix}0&1\\-1&0\end{matrix} \right), X_2=\left(\begin{matrix}0&1\\1&0\end{matrix} \right), X_3=\left(\begin{matrix}1&0\\0&-1\end{matrix} \right)\right\},$$ we have
	\begin{equation} \label{eqsl2} [X_1,X_2]=2X_3,\;[X_3,X_1]=2X_2\esp [X_3,X_2]=2X_1. \end{equation}
	The Killing form has signature $(+,+,-)$.
	\item The solvable Lie group $\sol=\left\{\left(\begin{matrix} e^x&0&y\\0&e^{-x}&z\\0&0&1\end{matrix}\right),x,y,z\in\R   \right\}$ whose Lie algebra is $\mathfrak{sol}=\left\{\left(\begin{matrix} x&0&y\\0&-x&z\\0&0&0\end{matrix}\right),x,y,z\in\R   \right\}$. In the basis
	\[ \B_0=\left\{X_1=\left(\begin{matrix} 1&0&0\\0&-1&0\\0&0&0\end{matrix}\right),
	X_2=\left(\begin{matrix} 0&0&1\\0&0&0\\0&0&0\end{matrix}\right),\;
	 X_3=\left(\begin{matrix} 0&0&0\\0&0&1\\0&0&0\end{matrix}\right)\right\}, \]we have
	\begin{equation} \label{eqsol} [X_1,X_2]=X_2; [X_1,X_3]=-X_3\esp[X_2,X_3]=0. \end{equation}The Killing form has signature $(+,0,0)$.
	\item The universal covering group $\wi{\mathrm{E}_0}(2)$ of the Lie group 
	\[ \mathrm{E}_0(2)=\left\{\left(\begin{matrix} \cos(\theta)&\sin(\theta)&x\\-\sin(\theta)&\cos(\theta)&y\\0&0&1\end{matrix}\right),\theta,x,y\in\R   \right\}. \]Its Lie algebra is
	\[ \mathrm{e}_0(2)=\left\{\left(\begin{matrix} 0&\theta&x\\-\theta&0&y\\0&0&0\end{matrix}\right),\theta,y,z\in\R   \right\}.  \]In the basis
	\[ \B_0=\left\{X_1=\left(\begin{matrix} 0&-1&0\\1&0&0\\0&0&0\end{matrix}\right),
	X_2=\left(\begin{matrix} 0&0&1\\0&0&0\\0&0&0\end{matrix}\right),\;
	 X_3=\left(\begin{matrix} 0&0&0\\0&0&1\\0&0&0\end{matrix}\right)\right\}, \]we have
	\begin{equation} \label{eqe2} [X_1,X_2]=X_3,\; [X_1,X_3]=-X_2\esp[X_2,X_3]=0. \end{equation}
	The Killing form has signature $(-,0,0)$.
	
\end{enumerate}

\begin{rem}\label{rem2} We can see from what above that two unimodular three dimensional non abelian Lie algebras are isomorphic if and only their Killing forms have the same signature.
	
\end{rem}

\section{Lorentzian left invariant metrics on three dimensional unimodular Lie groups}
\label{section3}
Through this section, for any $G$ among the five unimodular Lie groups described in Section \ref{section2}, $\B_0$ is the basis of its Lie algebra given also in Section \ref{section2}, $h_0$ is a left invariant Lorentzian metric on $G$ and $L:\g\too\g$ the associated endomorphism given in \eqref{br}. 

\subsection{ Lorentzian left invariant metrics on $\nil$}
We have three possibilities:

\begin{enumerate}
	\item[$(i)$] $L$ is of type $\mathrm{diag}(0,0,\sqrt{\la})$ with $\la\not=0$. In the basis $\B_1$ given in \eqref{bd}, we have 
	$[e_1,e_2]=-\sqrt{\la} e_3.$
	We consider the automorphism of Lie algebras $P:\mathfrak{n}\too \mathfrak{n}$    given by
	\[ P(X)=e_1;\; P(Y)=e_2\esp P(Z)=-\sqrt{\la} e_3, \] $\phi:G\too G$ the associated automorphism of Lie groups. The matrix of   $\phi^*(h_0)$ is given by
	\begin{equation}\label{n1} M(\phi^*(h_0),\B_0)=\left(\begin{matrix}
	1&0&0\\0&1&0\\0&0&-\la
	\end{matrix}  \right),\quad\la>0.  \end{equation}
	\item[$(ii)$] $L$ is of type $\mathrm{diag}(\sqrt{\la},0,0)$ with $\la\not=0$. In the basis $\B_1$ given in \eqref{bd}, we have 
	$[e_2,e_3]=\sqrt{\la} e_1.$
	We consider the automorphism of Lie algebras $P:\mathfrak{n}\too \mathfrak{n}$ given by
	\[ P(X)=e_2;\; P(Y)=e_3\esp P(Z)=\sqrt{\la} e_1, \]$\phi$ its associated automorphism of Lie groups. The matrix of  $\phi^*(h_0)$ in $\B_0$ is given by
	\begin{equation}\label{n2} M(\phi^*(h_0),\B_0)= \left(\begin{matrix}
	\la&0&0\\0&1&0\\0&0&-1
	\end{matrix}  \right),\;\la>0.  \end{equation}
	
	\item[$(iii)$] $L$ is of type $\{002\}$.  In the basis $\B_1$ given in \eqref{bab2}, we have 
	$$[e_1,e_2]=[e_3,e_1]=\frac12 (e_2+e_3).$$
	We consider the automorphism of Lie algebras $P:\mathfrak{n}\too \mathfrak{n}$ given by
	\[ P(X)=e_1;\; P(Y)=e_2\esp P(Z)=\frac12( e_2+e_3), \]$\phi$ its associated automorphism of Lie groups. The matrix of  $\phi^*(h_0)$ is given by 
	\begin{equation}\label{n0}
	M(\phi^*(h_0),\B_0)=\left(\begin{matrix}
	1&0&0\\0&1&\frac12\\0&\frac12&0
	\end{matrix}  \right).
	\end{equation}
\end{enumerate}
\begin{theo} \label{nil} Any  Lorentzian left invariant metric on $\nil$ is isometric to one of the three  metrics whose matrices in $\B_0$ are given by \eqref{n1}-\eqref{n0}.

\end{theo}

\subsection{ Lorentzian left invariant metric on $\mathrm{SU}(2)$}

We have one possibility, $L$ is of type $\mathrm{diag}\{a,b,c \}$
with $a>0$, $b>0$ and $c<0$. In the basis $\B_1$ given in \eqref{bd}, we have 
\[ [e_1,e_2]=-ce_3,\;[e_2,e_3]=ae_1\esp [e_3,e_1]=be_2. \]
We consider the automorphism $P:\mathfrak{su}(2)\too \mathfrak{su}(2)$ given by
\[ P(\sigma_x)=\frac{2}{\sqrt{-cb}}e_1=\sqrt{\mu_1}e_1;\; P(\sigma_y)=\frac{2}{\sqrt{-ca}}e_2=\sqrt{\mu_2}e_2\esp P(\sigma_z)=\frac{2}{\sqrt{ab}}e_3=\sqrt{\mu_3}e_3, \]$\phi$ its associated automorphism of Lie groups. The matrix of $\phi^*(h_0)$ in $\B_0$ is given 
\begin{equation}\label{su2}
M(\phi^*(h_0),\B_0)=\left(\begin{matrix}
\mu_1&0&0\\0&\mu_2&0\\0&0&-\mu_3
\end{matrix}  \right),\;\mu_1\geq\mu_2>0,\mu_3>0.
\end{equation}We can suppose that $\mu_1\geq\mu_2$ since  $L=\left(\begin{matrix}
0&1&0\\1&0&0\\0&0&-1
\end{matrix}  \right)$  is an automorphism of $\mathfrak{su}(2)$ and
$$L^tM(\phi^*(h_0),\B_0)L=\left(\begin{matrix}
\mu_2&0&0\\0&\mu_1&0\\0&0&-\mu_3
\end{matrix}  \right).$$

\begin{theo}\label{su} Any  Lorentzian left invariant metric on $\mathrm{SU}(2)$ is isometric to the metric whose matrix in $\B_0$ is given by \eqref{su2}.  
\end{theo}

\subsection{ Lorentzian left invariant metrics on $\wi{\mathrm{PSL}}(2,\R)$}

We have five possibilities:
\begin{enumerate}
	\item $L$ is of type $\mathrm{diag}\{a,b,c \}$
	with $a>0$, $b>0$ and $c>0$. In the basis $\B_1$ given in \eqref{bd}, we have
	\[ [e_1,e_2]=-ce_3,\;[e_2,e_3]=ae_1\esp [e_3,e_1]=be_2. \]
	We consider the automorphism $P:\mathrm{sl}(2,\R)\too \mathrm{sl}(2,\R)$ given by
	\[ P(X_1)=\frac{2}{\sqrt{ab}}e_3=\sqrt{\mu_1}e_3;\; P(X_2)=\frac{2}{\sqrt{ca}}e_2=\sqrt{\mu_2}e_2\esp P(X_3)=-\frac{2}{\sqrt{ab}}e_1=-\sqrt{\mu_3}e_1, \]$\phi$ its associated automorphism of Lie group. The matrix of  $\phi^*(h_0)$ in $\B_0$ is given by
	\begin{equation}\label{sl2d1}
	M(\phi^*(h_0),\B_0)=\left(\begin{matrix}
	-\mu_1&0&0\\0&\mu_2&0\\0&0&\mu_3
	\end{matrix}  \right),\; \quad \mu_1>0,\mu_2\geq\mu_3>0.
	\end{equation}We can suppose that $\mu_2\geq\mu_3$ since  $L=\left(\begin{matrix}
	1&0&0\\0&0&1\\0&-1&0
	\end{matrix}  \right)$  is an automorphism of $\mathrm{sl}(2,\R)$ and
	$$L^tM(\phi^*(h_0),\B_0)L=\left(\begin{matrix}
	-\mu_1&0&0\\0&\mu_3&0\\0&0&\mu_2
	\end{matrix}  \right).$$

	\ \item $L$ is of type $\mathrm{diag}\{a,b,c \}$
	with $a<0$, $b>0$ and $c<0$. 
	\[ [e_1,e_2]=-ce_3,\;[e_2,e_3]=ae_1\esp [e_3,e_1]=be_2. \]
	We consider  the automorphism $P:\mathrm{sl}(2,\R)\too \mathrm{sl}(2,\R)$ given by
	\[ P(X_1)=\frac{2}{\sqrt{-cb}}e_1=\sqrt{\mu_1}e_1;\; P(X_2)=\frac{2}{\sqrt{-ab}}e_3=\sqrt{\mu_2}e_3\esp P(X_3)=-\frac{2}{\sqrt{ac}}e_2=-\sqrt{\mu_3}e_2 \]$\phi$ its associated automorphism of Lie group. The matrix of $\phi^*(h_0)$ is given by
	\begin{equation}\label{sl2d2}
	M(\phi^*(h_0),\B_0)=\left(\begin{matrix}
	\mu_1&0&0\\0&-\mu_2&0\\0&0&\mu_3
	\end{matrix}  \right),\quad \mu_1>0,\mu_2>0,\mu_3>0.
	\end{equation}
	\item $L$ is of type $\{a^2z\bar{z}\}$ $a\not=0$. In the basis $\B_1$ given in \eqref{bazz}, we have
	\[ [e_1,e_2]=-\beta e_2+\alpha e_3,\;[e_2,e_3]=a^2e_1\esp [e_3,e_1]=\al e_2-\beta e_3,\;\be\not=0. \]We distinguish three cases:
	\begin{enumerate}
		\item  $\al>0$. We consider  the automorphism $P:\mathrm{sl}(2,\R)\too \mathrm{sl}(2,\R)$ given by{\small
			\[ P(X_1)=\frac{2}{a\sqrt{\al(\al^2+\be^2)}}(\be e_2+\al e_3);\; P(X_2)=\frac{2}{a\sqrt{\alpha}}e_2=\esp 
			P(X_3)=-\frac{2}{\sqrt{\alpha^2+\be^2}}e_1, \]}$\phi$ its associated automorphism of Lie group. The matrix of $\phi^*(h_0)$ is given by
		\begin{equation}\label{sl2azz+}
		M(\phi^*(h_0),\B_0)=\frac{4}{a^2\al\sqrt{\al^2+\be^2}}\left(\begin{matrix}
		\frac{\be^2-\al^2}{\sqrt{\al^2+\be^2}}&\be&0\\\be&\sqrt{\al^2+\be^2}&0\\
		0&0&\frac{a^2\al}{\sqrt{\al^2+\be^2}}
		\end{matrix}  \right),\quad \al>0, \be>0.
		\end{equation}
We can suppose that $\be>0$ since  $L=\left(\begin{matrix}
	1&0&0\\0&-1&1\\0&0&-1
	\end{matrix}  \right)$  is an automorphism of $\mathrm{sl}(2,\R)$ and
	$$L^tM(\phi^*(h_0),\B_0)L=\frac{4}{a^2\al\sqrt{\al^2+\be^2}}\left(\begin{matrix}
		\frac{\be^2-\al^2}{\sqrt{\al^2+\be^2}}&-\be&0\\-\be&\sqrt{\al^2+\be^2}&0\\
		0&0&\frac{a^2\al}{\sqrt{\al^2+\be^2}}
		\end{matrix}  \right).
		.$$ 		
		
		\item  $\al<0$. We consider  the automorphism $P:\mathrm{sl}(2,\R)\too \mathrm{sl}(2,\R)$ given by{\small
			\[P(X_1)=\frac{2}{a\sqrt{|\alpha|}}e_2,
			P(X_2)=-\frac{2}{\sqrt{\alpha^2+\be^2}}e_1\esp P(X_3)=\frac{2}{a\sqrt{|\al|(\al^2+\be^2)}}(-\be e_2-\al e_3) 
			 \]}$\phi$ its associated automorphism of Lie group. The matrix of $\phi^*(h_0)$ is given by
					\begin{equation}\label{sl2azz-}
		M(\phi^*(h_0),\B_0)=\frac{4}{a^2\al\sqrt{\al^2+\be^2}}\left(\begin{matrix}
		-\sqrt{\al^2+\be^2}&0&\be\\0&\frac{a^2\al}{\sqrt{\al^2+\be^2}}&0\\
		\be&0&\frac{\al^2-\be^2}{\sqrt{\al^2+\be^2}}
		\end{matrix}  \right),\;\al<0,\be>0.
		\end{equation}We can suppose $\be>0$ by using a same argument as in the precedent case.
		
		\item 	$\al=0$. We consider  the automorphism $P:\mathrm{sl}(2,\R)\too \mathrm{sl}(2,\R)$ given by{\small
			\[P(X_1)=\frac{1}{\be}e_2+\frac{2}{a^2}e_3,
			P(X_2)=\frac{2}{\be}e_1\esp P(X_3)=\frac{1}{\be}e_2-\frac{2}{a^2}e_3 
			 \]}$\phi$ its associated automorphism of Lie group. The matrix of $\phi^*(h_0)$  is given
		$$
		M(\phi^*(h_0),\B_0)=\frac{1}{a^4\be^2}\left(\begin{matrix}
		a^4-4\be^2&0&a^4+4\be^2\\0&4a^4&0\\
		a^4+4\be^2&0&a^4-4\be^2
		\end{matrix}  \right).
		$$By putting $u=a^4-4\be^2$ and $v=a^4+4\be^2$, we get
		\begin{equation}\label{sl2azz0}
		M(\phi^*(h_0),\B_0)=\frac{16}{v^2-u^2}\left(\begin{matrix}
		u&0&v\\0&2(u+v)&0\\
		v&0&u
		\end{matrix}  \right),\quad v>0,v>u.
		\end{equation}
		
	\end{enumerate}
	
	\item $L$ est of type $\{ab2\}$ with $a\not=0$ and $b\not=0$. In the basis $\B_1$ given in \eqref{bab2}, we have
	\[ [e_1,e_2]=\frac12 e_2+\left(\frac12-b\right)e_3, 
	[e_2,e_3]=ae_1\esp [e_3,e_1]=\left(b+\frac12\right)e_2+\frac12e_3. \]
	We consider  the automorphism $P:\mathrm{sl}(2,\R)\too \mathrm{sl}(2,\R)$ given by
	\begin{eqnarray*}
		P(X_1)&=&\left(\frac12+\frac1{4b}-\frac2{ab}\right)e_2+
		\left(-\frac12+\frac1{4b}-\frac2{ab}\right)e_3,\\
		P(X_2)&=&\left(-\frac12-\frac1{4b}-\frac2{ab}\right)e_2+
		\left(\frac12-\frac1{4b}-\frac2{ab}\right)e_3,\\
		P(X_3)&=&-\frac{2}{b}e_1
	\end{eqnarray*}
	$\phi$ its associated automorphism of Lie group. The matrix of $\phi^*(h_0)$ is given by 
	\begin{equation}\label{sl2ab2}
	M(\phi^*(h_0),\B_0)=\frac1{2ab}\left(\begin{matrix}
	a-8&-a&0\\-a&a+8&0\\
	0&0&\frac{8a}{b}
	\end{matrix}  \right),\quad a\not=0, b\not=0.
	\end{equation}
	\item $L$ est of type $\{a3\}$ with $a\not=0$.  In the basis $\B_1$ given in \eqref{ba3}, we have
	\begin{eqnarray*}
		[e_1,e_2]&=&\frac1{\sqrt{2}}e_2-ae_3,\;[e_2,e_3]=ae_1+\frac1{\sqrt{2}}e_2\esp \;[e_3,e_1]=\frac1{\sqrt{2}}e_1+ae_2+\frac1{\sqrt{2}}e_3.
	\end{eqnarray*}
	We consider  the automorphism $P:\mathrm{sl}(2,\R)\too \mathrm{sl}(2,\R)$ given by
	\begin{eqnarray*}
		P(X_1)&=&\frac1{a^2}\left(\sqrt{2}e_2-2ae_3\right),
		P(X_2)=\frac1{a^2\sqrt{2a^2+1}}\left(2ae_1+\sqrt{2}(2a^2+1)e_2\right),\\
		P(X_3)&=&\frac{2\sqrt{2}}{\sqrt{2a^2+1}}e_1
	\end{eqnarray*}
	$\phi$ its associated automorphism of Lie group. The matrix of $\phi^*(h_0)$ in $\B_0$ is given by 
	\begin{equation}\label{sl2a3}
	M(\phi^*(h_0),\B_0)=\frac2{a^4(1+2a^2)}\left(\begin{matrix}
	1-4a^4&(1+2a^2)^{\frac32}&0\\(1+2a^2)^{\frac32}&4a^4+6a^2+1&2a^3\sqrt{2}\\
	0&2a^3\sqrt{2}&4a^4
	\end{matrix}  \right). \quad a\not=0.
	\end{equation}

\end{enumerate}

\begin{theo}\label{sl2} Any Lorentzian left invariant metric on $\wi{\mathrm{PSL}}(2,\R)$ is isometric to one of the seven metrics whose matrices in $\B_0$ are given by
\eqref{sl2d1}-\eqref{sl2a3}.	
	
\end{theo}

\subsection{ Lorentzian left invariant metrics on $\sol$}

We have five possibilities:
\begin{enumerate}\item $L$ is of type $\mathrm{diag}(\alpha,\be,0)$ with $\alpha>0$, $\be<0$.   In the basis $\B_1$ given in \eqref{bd}, we have
	\begin{eqnarray*}
		[e_1,e_2]&=&0,\;[e_2,e_3]=\alpha e_1\esp \;[e_3,e_1]=\beta e_2,\quad \alpha>0,\beta<0.
	\end{eqnarray*}
	We consider  the automorphism $P:\mathfrak{sol}\too \mathfrak{sol}$ given by
	\begin{eqnarray*}
		P(X_1)&=&-\frac1{\sqrt{-\al\be}}e_3,
		P(X_2)=e_1-\frac{\be}{\sqrt{-\al\be}}e_2,
		P(X_3)=e_1+\frac{\be}{\sqrt{-\al\be}}e_2
	\end{eqnarray*}
	$\phi_0$ its associated automorphism of Lie group and we put $h_1=\phi_0^*(h_0)$. We have
	$$
	M(h_1,\B_0)=\frac1{\al}\left(\begin{matrix}
	\frac1\be&0&0\\0&\al-\be&\al+\be\\
	0&\al+\be&\al-\be
	\end{matrix}  \right).
	$$We can reduce this metric by considering the automorphism  $Q:\mathfrak{sol}\too \mathfrak{sol}$ given by
	\[ M(Q,\B_0)=\left( \begin{array}{ccc}
	-1&0&0\\0&0&\sqrt{\frac{\al}{\al-\be}}\\
	0&\sqrt{\frac{\al}{\al-\be}}&0
	\end{array}     \right). \]  Consider $\phi$  the automorphism of $\sol$ associated to $Q$. The matrix of  $\phi^*(h_1)$ is given by 
	$$
	M(\phi^*(h_1),\B_0)=M(Q,\B_0)^tM(h_1,\B_0)M(Q,\B_0)=\left(\begin{matrix}
	\frac1{\al\be}&0&0\\0&1&\frac{\al+\be}{\al-\be}\\
	0&\frac{\al+\be}{\al-\be}&1
	\end{matrix}  \right).
	$$
We put $u=\al+\be$ and $v=\al-\be$ and we get	
	\begin{equation}\label{sold1}
	M(\phi^*(h_1),\B_0)=\left(\begin{matrix}
	\frac4{u^2-v^2}&0&0\\0&1&\frac{u}{v}\\
	0&\frac{u}{v}&1
	\end{matrix}  \right), v>0, u<v.
	\end{equation}
	\item $L$ is of type $\mathrm{diag}(a,b,c)$ with $a>0$, $b=0$ and $c>0$.  In the basis $\B_1$ given in \eqref{bd}, we have
	\begin{eqnarray*}
		[e_1,e_2]&=&\alpha e_3,\;[e_2,e_3]=\beta e_1\esp \;[e_3,e_1]=0,\; \alpha<0,\beta>0.
	\end{eqnarray*}
	We consider  the automorphism $P:\mathfrak{sol}\too \mathfrak{sol}$ given by
	\begin{eqnarray*}
		P(X_1)&=&-\frac1{\sqrt{-\al\be}}e_2,
		P(X_2)=e_3-\frac{\be}{\sqrt{-\al\be}}e_1,
		P(X_3)=e_3+\frac{\be}{\sqrt{-\al\be}}e_1
	\end{eqnarray*}
	$\phi_0$ its associated automorphism of Lie group and we put $h_1=\phi_0^*(h_0)$. We have
	$$
	M(h_1,\B_0)=\frac1{\al}\left(\begin{matrix}
	-\frac1\be&0&0\\0&-\al-\be&-\al+\be\\
	0&-\al+\be&-\al-\be
	\end{matrix}  \right).
	$$We can reduce this metric by considering the automorphism  $Q:\mathfrak{sol}\too \mathfrak{sol}$ given by
	\[ M(Q,\B_0)=\left( \begin{array}{ccc}
	-1&0&0\\0&0&\sqrt{\frac{\al}{\al-\be}}\\
	0&\sqrt{\frac{\al}{\al-\be}}&0
	\end{array}     \right). \]  Consider $\phi$  the automorphism of $\sol$ associated to $Q$. The matrix of   $\phi^*(h_1)$ in $\B_0$ is given by
	$$
	M(\phi^*(h_1),\B_0)=M(Q,\B_0)^tM(h_1,\B_0)M(Q,\B_0)=\left(\begin{matrix}
	-\frac1{\al\be}&0&0\\0&\frac{\al+\be}{\be-\al}&-1\\
	0&-1&\frac{\al+\be}{\be-\al}
	\end{matrix}  \right).
	$$We put $u=\al+\be$ and $v=\be-\al$ and we get	
	\begin{equation}\label{sold2}
	M(\phi^*(h_1),\B_0)=\left(\begin{matrix}
	\frac4{v^2-u^2}&0&0\\0&\frac{u}{v}&-1\\
	0&-1&\frac{u}{v}&
	\end{matrix}  \right), v>0, u<v.
	\end{equation}

	\item $L$ is of type $\{0z\bar{z}\}$.  In the basis $\B_1$ given in \eqref{bazz}, we have
	\begin{eqnarray*}
		[e_1,e_2]&=&-\beta e_2-\alpha e_3,\;[e_2,e_3]=0\esp \;[e_3,e_1]=\al e_2-\be e_3,\quad \beta\not=0.
	\end{eqnarray*}
	We distinguish two cases:\begin{enumerate}
		\item $\al\not=0$. We consider  the automorphism $P:\mathfrak{sol}\too \mathfrak{sol}$ given by
		\begin{eqnarray*}
			P(X_1)&=&\frac1{\sqrt{\al^2+\be^2}}e_1,
			P(X_2)=e_3+\frac{\be-\sqrt{\al^2+\be^2}}{\al}e_2,
			P(X_3)=e_3+\frac{\be+\sqrt{\al^2+\be^2}}{\al}e_2
		\end{eqnarray*}
		$\phi_0$ its associated automorphism of Lie group and we put $h_1=\phi_0^*(h_0)$. We have
		$$
		M(h_1,\B_0)=\left(\begin{matrix}
		\frac1{\al^2+\be^2}&0&0\\0&-\frac{2\be(-\be+\sqrt{\al^2+\be^2})}{\al^2}&-2\\
		0&-2&\frac{2\be(\be+\sqrt{\al^2+\be^2})}{\al^2}
		\end{matrix}  \right).
		$$We can suppose that $\be>0$ by using the automorphism of $\mathfrak{sol}$  given by $L=\left( \begin{matrix} -1&0&0\\0&0&1\\ 0&1&0  \end{matrix} \right)$.
		
		We can reduce this metric by considering the automorphism   $Q:\mathfrak{sol}\too \mathfrak{sol}$ given by
		\[ M(Q,\B_0)=\left( \begin{array}{ccc}
		1&0&0\\0&-\frac{\sqrt{2}\sqrt{\beta(\beta+\sqrt{\beta^2+\al^2})}}{2\alpha}&0\\
		0&0&\frac{\al\sqrt{2}}{2\sqrt{\beta(\beta+\sqrt{\beta^2+\al^2})}}
		\end{array}     \right). \]  Consider $\phi$  the automorphism of $\sol$ associated to $Q$. The matrix of  $\phi^*(h_1)$ in $\B_0$ is given by
		\begin{equation}\label{sol0zz}
		M(\phi^*(h_1),\B_0)=\left(\begin{matrix}
		\frac1{u+v}&0&0\\0&-\frac{v}{u}&1\\
		0&1&1
		\end{matrix}  \right),\;u=\al^2>0, v=\be^2>0.
		\end{equation}
		
		\item $\al=0$. We consider  the automorphism $P:\mathfrak{sol}\too \mathfrak{sol}$ given by
		\begin{eqnarray*}
			P(X_1)&=&\frac1{\be}e_1,
			P(X_2)=e_3,
			P(X_3)=e_2
		\end{eqnarray*}
		$\phi$ its associated automorphism of Lie group. We have 
		\begin{equation}\label{sol0zz0}
		M(\phi^*(h_0),\B_0)=\left(\begin{matrix}
		\frac1{u}&0&0\\0&-1&0\\
		0&0&1
		\end{matrix}  \right),\;u=\be^2>0
		\end{equation}
	\end{enumerate}
	\item $L$ is of type $\{a02\}$ with $a<0$.  In the basis $\B_1$ given in \eqref{bab2}, we have
	\[ [e_1,e_2]=\frac12 e_2+\frac12e_3, 
	[e_2,e_3]=-a^2e_1\esp [e_3,e_1]=\frac12e_2+\frac12e_3. \]
	We consider  the automorphism $P:\mathfrak{sol}\too \mathfrak{sol}$ given by
	\begin{eqnarray*}
		P(X_1)&=&\frac{\sqrt{2}}{a}e_3,
		P(X_2)=a\sqrt{2} e_1+e_2+e_3,
		P(X_3)=-a\sqrt{2} e_1+e_2+e_3
	\end{eqnarray*}
	$\phi_0$ its associated automorphism of Lie group and we put $h_1=\phi_0^*(h_0)$. We have
	$$
	M(h_1,\B_0)=-\frac1{a^2}\left(\begin{matrix}
	2&a\sqrt{2}&a\sqrt{2}\\a\sqrt{2}&-2a^4&2a^4\\
	a\sqrt{2}&2a^4&-2a^4
	\end{matrix}  \right).
	$$We can reduce this metric by considering the automorphism  $Q:\mathfrak{sol}\too \mathfrak{sol}$ given by
	\[ M(Q,\B_0)=\left( \begin{array}{ccc}
	1&0&0\\-\frac{\sqrt{2}(3-2a^2)}{8a^2}&-\frac{\sqrt{2}}{2a}&0\\
	-\frac{\sqrt{2}(1+2a^2)}{8a^3}&0&\frac{\sqrt{2}}{2a}
	\end{array}     \right). \]  Consider $\phi$  the automorphism of $\sol$ associated to $Q$. The matrix of   $\phi^*(h_1)$ is given by
	\begin{equation}\label{sola02}
	M(\phi^*(h_1),\B_0)=\left(\begin{matrix}
	0&0&-\frac{2}{b}\\0&1&1\\
	-\frac{2}{b}&1&1
	\end{matrix}  \right),\quad b=a^2>0.
	\end{equation}
	
	\item $L$ is of type $\{0b2\}$ with $b\not=0$.  In the basis $\B_1$ given in \eqref{bab2}, we have
	\[ [e_1,e_2]=\frac12 e_2+\left(\frac12-b\right)e_3, 
	[e_2,e_3]=0\esp [e_3,e_1]=\left(b+\frac12\right)e_2+\frac12e_3. \]
	We consider  the automorphism $P:\mathfrak{sol}\too \mathfrak{sol}$ given by
	\begin{eqnarray*}
		P(X_1)&=&\frac{2}{b}e_1,
		P(X_2)=(2b+1) e_2+(1-2b)e_3,
		P(X_3)=e_2+e_3
	\end{eqnarray*}
	$\phi_0$ its associated automorphism of Lie group and we put $h_1=\phi_0^*(h_0)$. We have
	$$
	M(h_1,\B_0)=\left(\begin{matrix}
	\frac1{b^2}&0&0\\0&8b&4b\\
	0&4b&0
	\end{matrix}  \right).
	$$We can reduce this metric by considering the automorphism  $Q:\mathfrak{sol}\too \mathfrak{sol}$ given by
	\[ M(Q,\B_0)=\left( \begin{array}{ccc}
	1&0&0\\0&\frac{\sqrt{2}}{4b}&0\\
	0&0&\frac{\sqrt{2}}{2}
	\end{array}     \right). \]  Consider $\phi$  the automorphism of $\sol$ associated to $Q$. The matrix of   $\phi^*(h_1)$ is given by
	\begin{equation}\label{sol0b2}
	M(\phi^*(h_1),\B_0)=\left(\begin{matrix}
	\la^2&0&0\\0&\la&1\\
	0&1&0
	\end{matrix}  \right),\quad \la=\frac{1}{b}\not=0.
	\end{equation}

	\item $L$ est of type $\{03\}$. In the basis $\B_1$ given in \eqref{ba3}, we have
	\begin{eqnarray*}
		[e_1,e_2]&=&\frac1{\sqrt{2}}e_2,\;[e_2,e_3]=\frac1{\sqrt{2}}e_2\esp \;[e_3,e_1]=\frac1{\sqrt{2}}e_1+\frac1{\sqrt{2}}e_3.
	\end{eqnarray*}
	We consider  the automorphism $P:\mathfrak{sol}\too \mathfrak{sol}$ given by
	\begin{eqnarray*}
		P(X_1)&=&\sqrt{2}e_1,
		P(X_2)=e_2\esp
		P(X_3)=e_3
	\end{eqnarray*}
	$\phi_0$ its associated automorphism of Lie group and we put $h_1=\phi_0^*(h_0)$. We have
	$$
	M(h_1,\B_0)=\left(\begin{matrix}
	2&0&\sqrt{2}\\0&1&0\\\sqrt{2}&0&0
	\end{matrix}  \right).
	$$We can reduce this metric by considering the automorphism  $Q:\mathfrak{sol}\too \mathfrak{sol}$ given by
	\[ M(Q,\B_0)=\left( \begin{array}{ccc}
	1&0&0\\0&1&0\\
	-\frac{\sqrt{2}}{2}&0&\frac{\sqrt{2}}{2}
	\end{array}     \right). \]  Consider $\phi$  the automorphism of $\sol$ associated to $Q$. The matrix of  $\phi^*(h_1)$ is given by
	\begin{equation}\label{s0l03}
	M(\phi^*(h_1),\B_0)=\left(\begin{matrix}
	0&0&1\\0&1&0\\1&0&0
	\end{matrix}  \right).
	\end{equation}

\end{enumerate}
\begin{theo}\label{sol} Any Lorentzian left invariant metric on $\sol$ is isometric to one of the six metrics whose matrices in $\B_0$ are given by
	\eqref{sold1}-\eqref{s0l03}.	
	
\end{theo}

\subsection{ Lorentzian left invariant metrics on $\wi{\mathrm{E}_0}(2)$}
There are three possibilities: 
\begin{enumerate}
	\item $L$ is of type $\mathrm{diag}(a,b,c)$ with $a>0$, $b>0$ and $c=0$.  In the basis $\B_1$ given in \eqref{bd}, we have
	\begin{eqnarray*}
		[e_1,e_2]&=&0,\;[e_2,e_3]=\alpha e_1\esp \;[e_3,e_1]=\beta e_2,\quad \alpha>0,\beta>0.
	\end{eqnarray*}
	We consider  the automorphism $P:\mathrm{e}_0(2)\too \mathrm{e}_0(2)$ given by
	\begin{eqnarray*}
		P(X_1)&=&\frac1{\sqrt{\al\be}}e_3,
		P(X_2)=e_1-\frac{\be}{\sqrt{\al\be}}e_2,
		P(X_3)=e_1+\frac{\be}{\sqrt{\al\be}}e_2
	\end{eqnarray*}
	$\phi_0$ its associated automorphism of Lie group and we put $h_1=\phi_0^*(h_0)$. We have
	$$
	M(h_1,\B_0)=\frac1{\al}\left(\begin{matrix}
	-\frac1\be&0&0\\0&\al+\be&\al-\be\\
	0&\al-\be&\al+\be
	\end{matrix}  \right).
	$$ We can reduce this metric by considering the automorphism  $Q:\mathrm{e}_0(2)\too \mathrm{e}_0(2)$ given by
	\[ M(Q,\B_0)=\left( \begin{array}{ccc}
	1&0&0\\\frac{1}{2\sqrt{\al\be}}&\frac{\sqrt{\al\be}}{2}&-\frac{\sqrt{\al\be}}{2}\\
	\frac{1}{2\sqrt{\al\be}}&\frac{\sqrt{\al\be}}{2}&\frac{\sqrt{\al\be}}{2}
	\end{array}     \right). \]  Consider $\phi$  the automorphism of $\wi{\mathrm{E}_0}(2)$ associated to $Q$. The matrix of   $\phi^*(h_1)$ is given by
	\begin{equation}\label{e2d1}
	M(\phi^*(h_1),\B_0)=\left(\begin{matrix}
	0&1&0\\1&u&0\\
	0&0&v
	\end{matrix}  \right),\;u=\al\be>0,v=\be^2>0.
	\end{equation}

	\item $L$ is of type $\mathrm{diag}(a,b,c)$ with $a>0$, $b=0$ and $c<0$.  In the basis $\B_1$ given in \eqref{bd}, we have
	\begin{eqnarray*}
		[e_1,e_2]&=&\alpha e_3,\;[e_2,e_3]=\beta e_1\esp \;[e_3,e_1]=0,\; \alpha>0,\beta>0.
	\end{eqnarray*}
	We consider  the automorphism $P:\mathrm{e}_0(2)\too \mathrm{e}_0(2)$ given by
	\begin{eqnarray*}
		P(X_1)&=&\frac1{\sqrt{\al\be}}e_2,
		P(X_2)=e_3-\frac{\be}{\sqrt{\al\be}}e_1,
		P(X_3)=e_3+\frac{\be}{\sqrt{\al\be}}e_1
	\end{eqnarray*}
	$\phi_0$ its associated automorphism of Lie group and we put $h_1=\phi_0^*(h_0)$. We have
	$$
	M(h_1,\B_0)=\frac1{\al}\left(\begin{matrix}
	\frac1\be&0&0\\0&-\al+\be&-\al-\be\\
	0&-\al-\be&-\al+\be
	\end{matrix}  \right).
	$$We can reduce this metric by considering the automorphism  $Q:\mathrm{e}_0(2)\too \mathrm{e}_0(2)$ given by
	\[ M(Q,\B_0)=\left( \begin{array}{ccc}
	1&0&0\\\frac{1}{2\sqrt{\al\be}}&\frac{\sqrt{\al\be}}{2}&-\frac{\sqrt{\al\be}}{2}\\
	\frac{1}{2\sqrt{\al\be}}&\frac{\sqrt{\al\be}}{2}&\frac{\sqrt{\al\be}}{2}
	\end{array}     \right). \]  Consider $\phi$  the automorphism of $\wi{\mathrm{E}_0}(2)$ associated to $Q$. The matrix of  $\phi^*(h_1)$ is given by
	\begin{equation}\label{e2d2}
	M(\phi^*(h_1),\B_0)=\left(\begin{matrix}
	0&-1&0\\-1&-u&0\\
	0&0&v
	\end{matrix}  \right),\quad u=\al\be>0, v=\be^2>0.
	\end{equation}

	\item $L$ is of type $\{a02\}$ with $a>0$.  In the basis $\B_1$ given in \eqref{bab2}, we have
	\[ [e_1,e_2]=\frac12 e_2+\frac12e_3, 
	[e_2,e_3]=a^2e_1\esp [e_3,e_1]=\frac12e_2+\frac12e_3. \]
	We consider  the automorphism $P:\mathrm{e}_0(2)\too \mathrm{e}_0(2)$ given by
	\begin{eqnarray*}
		P(X_1)&=&\frac{\sqrt{2}}{a}e_3,
		P(X_2)=a\sqrt{2} e_1+e_2+e_3,
		P(X_3)=-a\sqrt{2} e_1+e_2+e_3
	\end{eqnarray*}
	$\phi_0$ its associated automorphism of Lie group and we put $h_1=\phi_0^*(h_0)$. We have
	$$
	M(h_1,\B_0)=-\frac1{a^2}\left(\begin{matrix}
	2&a\sqrt{2}&a\sqrt{2}\\a\sqrt{2}&-2a^4&2a^4\\
	a\sqrt{2}&2a^4&-2a^4
	\end{matrix}  \right).
	$$ We can reduce this metric by considering the automorphism  $Q:\mathrm{e}_0(2)\too \mathrm{e}_0(2)$ given by
	\[ M(Q,\B_0)=\left( \begin{array}{ccc}
	1&0&0\\-\frac{\sqrt{2}}{4a}&-\frac{\sqrt{2}a}{4}&-\frac{\sqrt{2}a}{4}\\
	-\frac{\sqrt{2}}{4a}&-\frac{\sqrt{2}a}{4}&-\frac{\sqrt{2}a}{4}
	\end{array}     \right). \]  Consider $\phi$  the automorphism of $\wi{\mathrm{E}_0}(2)$ associated to $Q$. The matrix of  $\phi^*(h_1)$ is given by
	\begin{equation}\label{e2a02}
	M(\phi^*(h_1),\B_0)=\left(\begin{matrix}
	0&1&0\\1&0&0\\
	0&0&u
	\end{matrix}  \right), \quad u=a^4>0.
	\end{equation}
	
\end{enumerate}

\begin{theo}  \label{e2}Any Lorentzian left invariant metric on $\wi{\mathrm{E}_0}(2)$ is isometric to one of the three metrics whose matrices in $\B_0$ are given by
	\eqref{e2d1}-\eqref{e2a02}.	
	
\end{theo}

\section{Curvature of Lorentzian left invariant metrics on unimodular three dimensional Lie groups}\label{section4}

In this section, we give for each metric whose matrix is given by one of the formulas \eqref{n1}-\eqref{e2a02} its Ricci tensor,  its signature and the scalar curvature.

\subsection{Curvature of Lorentzian left invariant metrics on $\nil$}

There are three classes of metrics on $\nil$ given by the formulas \eqref{n1}-\eqref{n0}.
Here are their Ricci curvature and scalar curvature.

\begin{pr} \begin{enumerate}
		\item The Ricci curvature and the scalar curvature of the metric \eqref{n1} are given by
		\[ M(\ric,\B_0)=\frac12\left(\begin{matrix}
		\la&0&0\\0&\la&0\\0&0&\la^2
		\end{matrix}  \right),\quad \mathfrak{s}=\frac12\la. \]In particular, $\ric>0$ and $\mathfrak{s}>0$.
	\item 	The Ricci curvature and the scalar curvature of the metric \eqref{n2} are given by
	\[ M(\ric,\B_0)=\frac1{2}\left(\begin{matrix}
	1&0&0\\0&\frac1{\la}&0\\0&0&\frac1{\la}
	\end{matrix}  \right),\quad \mathfrak{s}=\frac1{2\la}. \]In particular, $\ric>0$ and $\mathfrak{s}>0$.
	\item The metric \eqref{n0} is flat.
	\end{enumerate}
	
\end{pr}

\subsection{Curvature of Lorentzian left invariant metrics on $\mathrm{SU}(2)$}

There is one class of metrics on $\mathrm{SU}(2)$ given by the formula \eqref{su2}.
Here is its  Ricci curvature and scalar curvature.

\begin{pr} The Ricci curvature and the scalar curvature of the metric on $\mathrm{SU}(2)$ given by \eqref{su2}  are given by
	\begin{eqnarray*}
		M(\ric,\B_0)&=&\mathrm{diag}\left[-\frac{2(\mu_1-\mu_2-\mu_3)(\mu_1+\mu_2+\mu_3)}{
			\mu_2\mu_3},\frac{2(\mu_1+\mu_2+\mu_3)(\mu_1-\mu_2+\mu_3)}{
			\mu_1\mu_3},-\frac{2(\mu_1-\mu_2-\mu_3)(\mu_1-\mu_2+\mu_3)}{
			\mu_1\mu_2}\right],\\
		\mathfrak{s}&=&\frac{2((\sqrt{\mu_1}+\sqrt{\mu_2})^2+\mu_3)((\sqrt{\mu_1}-\sqrt{\mu_2})^2+\mu_3)}{
			\mu_1\mu_2\mu_3}>0.
	\end{eqnarray*}Moreover, the signature of $\ric$ is given by
	\[ \mathrm{sign}(\ric)=\left\{\begin{array}{ccc}
	(+,+,+)&\mbox{if}&\mu_1<\mu_2+\mu_3,\\
	(+,-,-)&\mbox{if}&\mu_1>\mu_2+\mu_3,\\
	(+,0,0)&\mbox{if}&\mu_1=\mu_2+\mu_3.
	\end{array} \right. \]
	
\end{pr}

\subsection{Curvature of Lorentzian left invariant metrics on $\wi{\mathrm{PSL}}(2,\R)$}
There are seven classes of metrics on $\wi{\mathrm{PSL}}(2,\R)$ given by the formulas \eqref{sl2d1}-\eqref{sl2a3}.
Here are their Ricci curvature and scalar curvature.

\begin{pr}
	The Ricci curvature and the scalar curvature of the metric \eqref{sl2d1} on $\wi{\mathrm{PSL}}(2,\R)$  are given by
	\begin{eqnarray*}
		M(\ric,\B_0)&=&\mathrm{diag}\left[
		\frac{2(\mu_1^2-(\mu_2-\mu_3)^2)}{
			\mu_2\mu_3},-\frac{2(\mu_2^2-(\mu_1-\mu_3)^2)}{
			\mu_1\mu_3},-\frac{2(\mu_3^2-(\mu_2-\mu_1)^2)}{
			\mu_1\mu_2}\right],\\
		\mathfrak{s}&=&\frac{2[(\sqrt{\mu_1}+\sqrt{\mu_2})^2-\mu_3]
			[(\sqrt{\mu_1}-\sqrt{\mu_2})^2-\mu_3]}{
			\mu_1\mu_2\mu_3}.
	\end{eqnarray*}When $\mu_1=\mu_2=\mu_3=\mu$ then $\ric=-\frac{2}\mu h$ and in fact the metric has constant sectional curvature $-\frac1{\mu}$. The possible signatures of the Ricci curvature are
	$(+,+,+),(+,-,-),(+,0,0),(-,0,0).$
\end{pr} 

\begin{pr} The Ricci curvature and the scalar curvature of the metric \eqref{sl2d2} on $\wi{\mathrm{PSL}}(2,\R)$  are given by
\begin{eqnarray*}
	M(\ric,\B_0)&=&\mathrm{diag}\left[-\frac{2(\mu_1-\mu_2-\mu_3)(\mu_1+\mu_2+\mu_3)}{
		\mu_2\mu_3},-\frac{2(\mu_1+\mu_2+\mu_3)(\mu_1-\mu_2+\mu_3)}{
		\mu_1\mu_3},\frac{2(\mu_1-\mu_2-\mu_3)(\mu_1-\mu_2+\mu_3)}{
		\mu_1\mu_2}\right],\\
	\mathfrak{s}&=&\frac{2((\sqrt{\mu_1}+\sqrt{\mu_2})^2+\mu_3)((\sqrt{\mu_1}-
		\sqrt{\mu_2})^2+\mu_3)}{
		\mu_1\mu_2\mu_3}>0.
\end{eqnarray*}	
Moreover, the signature of $\ric$ is given by
\[ \mathrm{sign}(\ric)=\left\{\begin{array}{ccc}
(+,+,+)&\mbox{if}&\mu_1<\mu_2-\mu_3<\mu_2+\mu_3,\\
(+,-,-)&\mbox{if}&\mu_1>\mu_2-\mu_3,\\
(-,0,0)&\mbox{if}&\mu_1=\mu_2-\mu_3\;\mbox{or}\;\mu_1=\mu_2+\mu_3.
\end{array} \right. \]

\end{pr}

\begin{pr} The Ricci curvature and the scalar curvature of the metric \eqref{sl2azz+} are given by
	\begin{eqnarray*}
		M(\ric,\B_0)&=&\left(\begin{matrix}
			\frac{2(a^2-2\al)(a^2(\be^2-\al^2)+4\al\be^2)}{a^2\al(\al^2+\be^2)}
			&\frac{2(a^4-4\al^2)\be}{a^2\al\sqrt{\al^2+\be^2}}&0\\
			\frac{2(a^4-4\al^2)\be}{a^2\al\sqrt{\al^2+\be^2}}&
			\frac{2(a^2-2\al)}{\al}&0\\
			0&0&-\frac{2(a^4-4\be^2)}{{\al^2+\be^2}}
		\end{matrix}  \right),\\
		\mathfrak{s}&=&\frac12a^4-2a^2\al-2\be^2.
	\end{eqnarray*}The signature of $\ric$ is $(+,-,-)$ if $a^2\not=2\al$ and $(-,0,0)$ if $a^2=2\al$. The operator of Ricci is of type $\{az\bar{z}\}$.
	
\end{pr}

\begin{pr}The Ricci curvature and the scalar curvature of the metric \eqref{sl2azz-}  are given by
	\begin{eqnarray*}
		M(\ric,\B_0)&=&\left(\begin{matrix}
			\frac{2(2\al-a^2)}{\al}
			&0&\frac{2(a^4-4\al^2)\be}{a^2\al\sqrt{\al^2+\be^2}}\\
			0&
			-\frac{2(a^4+4\be^2)}{\al^2+\be^2}&0\\
			\frac{2(a^4-4\al^2)\be}{a^2\al\sqrt{\al^2+\be^2}}&0&-\frac{2(a^2-2\al)(a^2(\be^2-\al^2)+4\al\be^2)}{a^2\al(\al^2+\be^2)}
		\end{matrix}  \right),\\
		\mathfrak{s}&=&\frac12a^4+2a^2\al-2\be^2.
	\end{eqnarray*}The signature of $\ric$ is $(+,-,-)$. The operator of Ricci is of type $\{az\bar{z}\}$.
	
\end{pr}

\begin{pr}The Ricci curvature and the scalar curvature of the metric \eqref{sl2azz0}  are given by
	\begin{eqnarray*}
		M(\ric,\B_0)&=&\left(\begin{matrix}
			\frac{4v}{v-u}
			&0&\frac{4v}{v-u}\\
			0&
			\frac{16v}{u-v}&0\\
			\frac{4v}{v-u}&0&\frac{4(v-2u)}{u-v}
		\end{matrix}  \right),\quad
		\mathfrak{s}=\frac{u}2.
	\end{eqnarray*}The signature of $\ric$ is $(+,-,-)$.
	
\end{pr}

\begin{pr} The Ricci curvature and the scalar curvature of the metric \eqref{sl2ab2} are given by
	\begin{eqnarray*}
		M(\ric,\B_0)&=&\frac1{4b}\left(\begin{matrix}
			(a+2b-8)(a-2b)&4b^2-a^2&0\\4b^2-a^2&(a+2b+8)(a-2b)&0\\
			0&0&-\frac{8a^2}{b}
		\end{matrix}  \right),\\
		\mathfrak{s}&=&\frac12a(a-4b).
	\end{eqnarray*}The Ricci curvature has signature $(+,-,-)$ if $a\not=2b$ and $(-,0,0)$ if $a=2b$.

\end{pr}

\begin{pr}
	The Ricci curvature and the scalar curvature of the metric \eqref{sl2a3} are given by
	\begin{eqnarray*}
		M(\ric,\B_0)&=&\left(\begin{array}{ccc}	\frac{2{a}^{2}
				-9}{{a}^{2}}&{\frac {-6\,{a}^{2}-9}{{a}^{2} \sqrt{2\,{a}^{2}+1}}}&-\,{\frac { 6\sqrt{2}}{ \sqrt{2\,{a}^{2}+1}a}}\\
			{\frac {-6\,{a}^{2}-9}{{a}^{2} \sqrt{2\,{a}^{2}+1}
					\mbox{}}}&{\frac {-4\,{a}^{4}-14\,{a}^{2}-9}{2\,{a}^{4}+{a}^{2}}}&-{\frac { \left( 6\,{a}^{2}+6 \right)  \sqrt{2}}{2\,{a}^{3}+a}}\\
			-\,{\frac { 6\sqrt{2}}{ \sqrt{2\,{a}^{2}+1}
					\mbox{}a}}&-{\frac { \left( 6\,{a}^{2}+6 \right)  \sqrt{2}}{2\,{a}^{3}+a}}&{\frac {-4\,{a}^{2}-8}{2\,{a}^{2}+1}}\end{array}
		\right),\quad
		\mathfrak{s}=-\frac32a^2.
	\end{eqnarray*}The operator $\Ric$ is of type $\{a3\}$. From the relations 
	$\det(M(\ric,\B_0))=8$ and $\tr(M(\ric,\B_0))=-2\frac{(a^2+9)}{a^2}$, we deduce that the signature of $\ric$ is $(+,-,-)$.
\end{pr}

\subsection{Curvature of Lorentzian left invariant metrics on $\sol$}

There are six classes of metrics on $\sol$ given by the formulas \eqref{sold1}-\eqref{s0l03}.
Here are their Ricci curvature and scalar curvature.

\begin{pr}
	The Ricci curvature and the scalar curvature of the metric \eqref{sold1} are given by
	\[M(\ric,\B_0)=\left(\begin{array}{ccc}\frac{2v^2}{u^2-v^2}&0&0\\
	0&-\frac12u^2&-\frac12uv\\0&-\frac12uv&-\frac12u^2   \end{array}\right),\quad
	\mathfrak{s}=\frac12v^2.  \]The signature of $\ric$ is given by
	\[\mathrm{sign}(\ric)= \left\{ \begin{array}{ccc}
	(-,0,0)&\mbox{if}&u=0,\\ 
	(-,-,0)&\mbox{if}&u=-v,\\
	(-,-,+)&\mbox{if}&u>0,\;\mbox{or}\; u<0,u>-v,\\
	(-,-,-)&\mbox{if}&u<0,u>-v.\end{array}       \right. \]
	
	\end{pr}
	\begin{pr}
		The Ricci curvature and the scalar curvature of the metric \eqref{sold2} are given by
		\[M(\ric,\B_0)=\left(\begin{array}{ccc}\frac{2u^2}{v^2-u^2}&0&0\\
		0&-\frac12uv&\frac12u^2\\0&\frac12u^2&-\frac12uv   \end{array}\right),\quad
		\mathfrak{s}=\frac12u^2.  \]This metric is flat if $u=0$. For $u\not=0$,
		the signature of $\ric$ is given by
		\[\mathrm{sign}(\ric)= \left\{ \begin{array}{ccc}
		(+,-,-)&\mbox{if}&u>0,\\ 
		(+,+,+)&\mbox{if}&u<0,u>-v,\\
		(+,+,-)&\mbox{if}&u<0,u<-v.\end{array}       \right. \]
		
	\end{pr}
\begin{pr}
	The Ricci curvature and the scalar curvature of the metric \eqref{sol0zz} are given by
	\[M(\ric,\B_0)=\left(\begin{array}{ccc}-\frac{2v}{v+u}&0&0\\
	0&2v&2v\\0&2v&-2u   \end{array}\right),\quad
	\mathfrak{s}=-2v.  \]The Ricci curvature has signature $(+,-,-)$ and the Ricci operator is of type $\{a,z\bar{z}\}$.
\end{pr}

\begin{pr}
	The Ricci curvature and the scalar curvature of the metric \eqref{sol0zz0} are given by
	\[M(\ric,\B_0)=\left(\begin{array}{ccc}-2&0&0\\
	0&0&0\\0&0&0   \end{array}\right),\quad
	\mathfrak{s}=-2u.  \]
\end{pr}

\begin{pr}
	The Ricci curvature and the scalar curvature of the metric \eqref{sola02} are given by
	\[M(\ric,\B_0)=\left(\begin{array}{ccc}-2&b&0\\
	b&-\frac12b^2&-\frac12b^2\\0&-\frac12b^2&-\frac12b^2   \end{array}\right),\quad
	\mathfrak{s}=\frac12b^2.  \]
	The Ricci operator is of type $\{ab2\}$. From the relations 
	$\det(M(\ric,\B_0))=\frac{b^4}2$ and $\tr(M(\ric,\B_0))=-2-b^2$, we deduce that the signature of $\ric$ is $(+,-,-)$.
\end{pr}
\begin{pr}
	The Ricci curvature and the scalar curvature of the metric \eqref{sol0b2} are given by
	\[M(\ric,\B_0)=\left(\begin{array}{ccc}0&0&0\\
	0&-\frac{2}\la&0\\0&0&0   \end{array}\right),\quad
	\mathfrak{s}=0.  \]We have $\Ric^2=0$. 
\end{pr}

\begin{pr}
	The Ricci curvature and the scalar curvature of the metric \eqref{s0l03} are given by
	\[M(\ric,\B_0)=\left(\begin{array}{ccc}-2&0&0\\
	0&0&0\\0&0&0   \end{array}\right),\quad
	\mathfrak{s}=0.  \]We have $\Ric^2=0$ and this metric is semi-symmetric.
\end{pr}

\subsection{Curvature of Lorentzian left invariant metrics on $\wi{\mathrm{E}_0}(2)$}
There are seven classes of metrics on $\wi{\mathrm{E}_0}(2)$ given by the formulas \eqref{e2d1}-\eqref{e2a02}.
Here are their Ricci curvature and scalar curvature.

\begin{pr}
	The Ricci curvature and the scalar curvature of the metric \eqref{e2d1} are given by
	\[M(\ric,\B_0)=\left(\begin{array}{ccc}\frac{v-u}v&\frac{v^2-u^2}{2v}&0\\
	\frac{v^2-u^2}{2v}&\frac{u(v^2-u^2)}{2v}&0\\0&0&\frac{u^2-v^2}{2}   \end{array}\right),\quad
	\mathfrak{s}=\frac{(u-v)^2}{2v}.  \]This metric is flat when $u=v$,
	\[ \mathrm{sign}(\ric)=(+,+,-)\;\mbox{if}\; u<v\esp 
	\mathrm{sign}(\ric)=(+,-,-)\;\mbox{if}\; u>v. \]
\end{pr}

\begin{pr}
	The Ricci curvature and the scalar curvature of the metric \eqref{e2d2} are given by
	\[M(\ric,\B_0)=\left(\begin{array}{ccc}\frac{v+u}v&\frac{u^2-v^2}{2v}&0\\
	\frac{u^2-v^2}{2v}&\frac{u(u^2-v^2)}{2v}&0\\0&0&\frac{u^2-v^2}{2}   \end{array}\right),\quad
	\mathfrak{s}=\frac{(u+v)^2}{2v}.  \]The signature of $\ric$ is given by
	\[ \mathrm{sign}(\ric)=(+,-,-)\;\mbox{if}\; u<v, 
	\mathrm{sign}(\ric)=(+,+,+)\;\mbox{if}\; u>v\esp\mathrm{sign}(\ric)=(+,0,0)
	\;\mbox{if}\; u=v. \]
\end{pr}

\begin{pr}
	The Ricci curvature and the scalar curvature of the metric \eqref{e2a02} are given by
	\[M(\ric,\B_0)=\left(\begin{array}{ccc}1&\frac{u}2&0\\
	\frac{u}2&0&0\\0&0&-\frac{u}2   \end{array}\right),\quad
	\mathfrak{s}=\frac{u}2.  \]The operator of Ricci is of type $\{ab2\}$ and the signature of $\ric$ is $(+,-,-)$.
\end{pr}

The following Table gives the possible signatures of Ricci curvature of Lorentzian left invariant metrics on unimodular three dimensional Lie groups and the metrics realizing these signatures.

\begin{center}
\begin{tabular}{|c|l|l|}
	\hline
	Signature of Ricci&Metrics realizing this signature&Remarks\\
	 curvature&&\\
	\hline
	$(0,0,0)$&\eqref{n0},[\eqref{sold2},$u=0$],[\eqref{e2d1},$u=v$]&These metrics are flat\\
	\hline
	$(+,+,+)$&\eqref{n1},\eqref{n2},\eqref{su2},\eqref{sl2d1},\eqref{sl2d2},\eqref{sold2},\eqref{e2d2}&\\
	\hline
	$(-,-,-)$&\eqref{sold1}&\\
	\hline
	$(+,-,-)$&\eqref{su2},\eqref{sl2d1},\eqref{sl2d2},\eqref{sl2azz+},\eqref{sl2azz-},
	\eqref{sl2azz0},\eqref{sl2ab2},\eqref{sl2a3},\eqref{sold1},\eqref{sold2}&The metric \eqref{sl2d1} has negative constant sectional \\
	&\eqref{sol0zz},\eqref{sola02},\eqref{e2d1},\eqref{e2d2},\eqref{e2a02}& curvature for $\mu_1=\mu_2=\mu_3$,\\
	&& The metric \eqref{sl2a3} is a shrinking Ricci soliton\\
	\hline
	$(+,+,-)$&\eqref{sold2},\eqref{e2d1}&\\
	\hline
	$(-,-,0)$&\eqref{sold1}&\\
	\hline
	$(+,0,0)$&\eqref{su2},\eqref{sl2d1},\eqref{sol0b2}&\\
	\hline
	$(-,0,0)$&\eqref{sl2d1},\eqref{sl2d2},\eqref{sl2azz+},\eqref{sl2ab2},\eqref{sold1},\eqref{sol0zz0},\eqref{sol0b2},\eqref{s0l03}&The metrics \eqref{sol0b2} and  \eqref{s0l03} are steady Ricci soliton   \\
	&&and semi-symmetric not locally symmetric,\\
	&&$[$\eqref{sl2ab2},$a=b\not=0$$]$ is a shrinking Ricci soliton\\
	\hline
\end{tabular}\end{center}

The following proposition gives the type of the Ricci operators.

\begin{pr}\label{pr}\begin{enumerate}
		\item The metrics \eqref{sl2azz+}, \eqref{sl2azz-} and \eqref{sol0zz} have their Ricci operators of type $\{az\bar{z}\}.$
		\item The metric \eqref{sl2a3} has its Ricci operator of type $\{a3\}.$
		\item The metrics \eqref{sol0b2} and \eqref{s0l03} have their Ricci operators of type $\{002\}$ and satisfies $\Ric^2=0$.
		\item The metrics \eqref{sola02} and \eqref{e2a02} have their Ricci operators of type $\{ab2\}.$
		\item All the others have their Ricci operators diagonalizable.
	\end{enumerate} 

\end{pr}

Lorentzian left invariant metrics on unimodular three dimensional Lie groups which are of constant curvature, Einstein, locally symmetric, semi-symmetric not locally symmetric or Ricci soliton have been determined in \cite{boucetta, bro, calvaruso1} by giving their Lie algebras as in \eqref{bd}-\eqref{ba3}. Our study permits a precise description of these metrics. We give  their corresponding metrics  in our list.
\begin{theo}[\cite{calvaruso1}]\label{main1} Let $h$ be a Lorentzian left invariant metric on a unimodular three dimensional Lie group. Then the following assertions are equivalent:
	\begin{enumerate}
		\item The metric $h$ is locally symmetric.
		\item The metric $h$ is Einstein.
		\item The metric $h$ has constant sectional curvature.
	\end{enumerate}
	Moreover, a metric satisfying one of these assertions  is either flat and is isometric to the metric \eqref{n0}, $[$\eqref{sold2},  $u=0$$]$ or $[$\eqref{e2d1},
	 $u=v$$]$; or it has a negative constant sectional curvature and is isometric to $[$\eqref{sl2d1},  $\mu_1=\mu_2=\mu_3$$]$.
\end{theo}

\begin{theo}[\cite{boucetta, calvaruso1}]\label{main2} The metrics \eqref{sol0b2} and \eqref{s0l03} are the unique, up to an automorphism, Lorentzian semi-symmetric not locally symmetric left invariant metrics on a unimodular three dimensional Lie group.
	
\end{theo}

\begin{theo}[\cite{bro}]\label{main3}\begin{enumerate}
		\item The metric \eqref{sol0b2} satisfies the relation
		\[ \mathrm{L}_Xh+\ric(h)=0,\quad X=-\frac1{\la^2}X_1 \]and hence it is a steady Ricci soliton.
		\item The metric \eqref{s0l03} satisfies the relation
		\[ \mathrm{L}_Xh+\ric(h)=0,\quad X=-X_3 \]and hence it is a steady Ricci soliton.
		\item The metric $[$\eqref{sl2ab2},$a=b\not=0$$]$ satisfies the relation
		\[ \mathrm{L}_Xh+\ric(h)=-\frac12b^2h,\quad X=\frac{b^2}{4}X_3 \]and hence it is a shrinking Ricci soliton.
		\item The metric \eqref{sl2a3} satisfies the relation
		\[ \mathrm{L}_Xh+\ric(h)=-\frac12a^2h,\quad X=\frac{a}{2\sqrt{2}}X_1
		+\frac{(2a^3-a)\sqrt{2}}{4\sqrt{2a^2+1}}X_2-
		\frac{a^2}{\sqrt{2a^2+1}}X3 \]and hence it is a shrinking Ricci soliton.
	\end{enumerate}Moreover, these metrics are the only Lorentzian left invariant Ricci solitons, up to automorphism, on three dimensional unimodular Lie groups.
	
\end{theo}

\end{document}